\documentclass[12pt,oneside]{amsart}

\usepackage{amsmath}
\usepackage{amsfonts}
\usepackage{amscd}
\usepackage{amsthm}
\usepackage{amssymb}
\usepackage{amsmath}
\usepackage{epsfig}
\usepackage{graphicx}
\usepackage[all]{xy}

\newcommand{\tM}{\tilde M}

\newcommand{\mD}{\mathbb D}
\newcommand{\mZ}{\mathbb Z}
\newcommand{\mN}{\mathbb N}

\newcommand{\mR}{\mathbb R}

\newcommand{\mU}{\mathcal U}
\newcommand{\om}{\omega}

\newcommand{\comment}[1]{}

\newcounter{fact}
\title{Large scale detection of half-flats in CAT(0) spaces.}
\author{Stefano Francaviglia}
\address{Departimento di Matematica, Universit\`a di Bologna, 
Piazza di Porta S. Donato 5, 40126, Bologna, Italy.}
\email{francavi@dm.unibo.it}
\author{Jean-Fran\c{c}ois Lafont}
\address{Department of Mathematics,
The Ohio State University, Columbus, OH 43210}
\email{jlafont@math.ohio-state.edu}

\theoremstyle{proposition}
\newtheorem{Lem}{Lemma}[section]

\newtheorem*{Def}{Definition}

\theoremstyle{plain}

\newtheorem{Cor}[Lem]{Corollary}
\newtheorem*{ThmA}{Theorem A}
\newtheorem*{ThmB}{Theorem B}

\newtheorem*{ThmC}{Theorem C}

\theoremstyle{remark}

\newtheorem*{Prf}{Proof}

\begin{document}

\begin{abstract}
Let $M$ be a complete locally compact CAT(0)-space, and $X$ an ultralimit of $M$.
For $\gamma \subset M$ a $k$-dimensional flat, let $\gamma_\omega$ be the
$k$-dimensional flat in $X$ obtained as an ultralimit of $\gamma$. 
In this paper, we identify various conditions
on $\gamma_\omega$ that are sufficient to ensure that $\gamma$ bounds a 
$(k+1)$-dimensional half-flat.

As applications we obtain (1) constraints on the behavior of quasi-isometries
between locally compact CAT(0)-spaces,
(2) constraints on the possible non-positively curved Riemannian metrics supported
by certain manifolds, and (3) a correspondence between metric splittings
of a complete, simply connected non-positively curved Riemannian manifolds, and metric 
splittings of its asymptotic cones.  Furthermore, combining our results with the 
Ballmann, Burns-Spatzier rigidity theorem and the classic Mostow rigidity,
we also obtain (4) a new proof of Gromov's rigidity theorem
for higher rank locally symmetric spaces.  
\end{abstract}

\maketitle

\section{Introduction.}

A {\it $k$-flat} in a CAT(0)-space $X$ is defined to be an isometrically embedded 
copy of the standard $\mR ^k$, $k\geq 1$. In the case where $k=1$, a $k$-flat is 
just a geodesic in $X$. 
By a $k$-dimensional {\it half-flat}, $k\geq 1$, in a CAT(0) space, 
we mean an isometric copy of $\mR ^{k-1}\times \mR ^+$ (where $\mR^+ = [0, \infty)$ 
is the usual half line). For example, when
$k=1$, a half-flat is just a geodesic ray in $X$. In the study of CAT(0)-spaces, a key role is 
played by the presence of flats and half-flats of {\it higher rank}, i.e. satisfying $k\geq 2$. In the 
present paper, our goal is to identify some coarse geometric conditions which
are sufficient to ensure the existence of half-flats in a CAT(0)-space $X$. We provide
three results towards this goal, as well as an example showing that our results are 
close to optimal. 

Before stating our main results, let us recall that an {\it asymptotic
cone} of a metric space $X$ is a new metric space, which encodes the large-scale
geometry of $X$, when viewed at an increasing sequence of scales. A precise
definition, along with some basic properties of asymptotic cones, is provided in
our Section 2. For a $k$-flat $\gamma_\omega$ inside the asymptotic cone $X$ of a 
CAT(0)-space, we introduce 
(see Section 3) the notion of a {\it flattening sequence} of maps for $\gamma_\omega$. These 
are a sequence of maps from a $k$-disk $\mD ^k$ into $X$, whose images are
getting further and further away from $\gamma_\omega$, and whose projections onto 
$\gamma_\omega$ satisfy certain technical conditions. The main point of such flattening 
sequences is that we can prove:

\begin{ThmA}[Flattening sequences $\Rightarrow$ half flats]
  Let $M$ be a locally compact CAT$(0)$-space and let $X$ be an
  asymptotic cone of $M$. Let $\gamma$ be a $k$-flat of $M$ (possibly
  a geodesic) and let $\gamma_\omega$ be its limit in $X$. Suppose
  that there exists a flattening sequence of maps for $\gamma _\omega$. Then
  the original $k$-flat $\gamma$ bounds a $(k+1)$-half-flat in $M$.
\end{ThmA}

The reader will readily see that, in the special case where $\gamma_\omega$
itself bounds a half-flat, it is very easy to construct a flattening sequence. So
an immediate consequence of Theorem A is:

\begin{ThmB}[Half-ultraflat $\Rightarrow$ half-flat] 
  Let $M$ be a locally compact CAT$(0)$-space and let $X$ be an
  asymptotic cone of $M$. Let $\gamma$ be a $k$-flat of $M$ (possibly
  a geodesic) and let $\gamma_\omega$ be its limit in $X$.
  If $\gamma_\omega$ bounds a half-flat, then $\gamma$ itself must
  bound a half-flat.
\end{ThmB}

\noindent We also provide an example showing that, in the context of locally compact
CAT(0)-spaces, the analogue of Theorem B with ``half-flats'' replaced by ``flats'' is {\it false}.

In Section 4, we weaken the hypothesis of Theorem B, by replacing a
$(k+1)$-flat in the ultralimit by a {\it bi-Lipschitzly embedded} $(k+1)$-flat. We
compensate for this by requiring the original flat to satisfy some mild periodicity
requirement, and  establish:

\begin{ThmC}[Bilipschitz half-ultraflat $+$ periodicity $\Rightarrow$ half-flat] 
  Let $M$ be a locally compact CAT$(0)$-space and let $X$ be an
  asymptotic cone of $M$. Let $\gamma$ be a $k$-flat of $M$ (possibly
  $k=1$) and let $\gamma_\omega$ its limit in $X$. Suppose that there
  exists $G<Isom(M)$ that acts co-compactly on $\gamma$. 
  
If there is a bi-Lipschitz embedding 
$\Phi:\mathbb R^k\times \mathbb R^+\to X$, 
whose restriction to $\mathbb R^k\times \{0\}$ is a homeomorphism onto
$\gamma_\omega$, then $\gamma$ bounds a $(k+1)$-half-flat in $M$. 
\end{ThmC}

In Section 5, we provide various geometrical 
applications of our main results.  These include:
\begin{itemize}
\item constraints on the possible quasi-isometries between certain
locally CAT(0)-spaces.
\item restrictions on the possible locally CAT(0)-metrics
that are supported by certain locally CAT(0)-spaces.
\item a proof that splittings of simply connected non-positively curved Riemannian 
manifolds correspond exactly with metric splittings of their asymptotic cones.
\item a new proof of Gromov's rigidity theorem \cite{BGS}: a closed higher rank
locally symmetric space supports a unique metric of non-positive curvature
(up to homothety).
\end{itemize}

\vskip 10pt

In January 2008, the authors posted a preliminary version \cite{FrLa} of this work on the arXiv,
which contained special cases of Thereoms A, B, C, under the additional hypothesis
that the flats be 2-dimensional, and the ambient space $M$ was a Riemannian 
manifold of non-positive sectional curvature (rather than a CAT(0)-space). 
Shortly thereafter, Misha Kapovich was kind
enough to inform the authors of his paper with B. Leeb \cite{KaL}, in which
(amongst other things) they proved a version of Theorem C in the special case of 
2-dimensional flats, and where the ambient space $M$ was an arbitrary locally compact 
CAT(0)-space. Their
paper provided the motivation for us to write the present paper, which includes a generalization to 
higher dimensional flats of the result in \cite[Prop. 3.3]{KaL}.

\vskip 10pt

\centerline{\bf Acknowledgements}

\vskip 10pt

The author's would like to thank V. Guirardel, J. Heinonen, T. Januszkiewicz,
B. Kleiner, R. Spatzier, and S. Wenger for their helpful comments. We are 
indebted to M. Kapovich for pointing out the existence of his joint paper with
B. Leeb \cite{KaL}, which put us on the track of the results contained in this paper.

This work was partly carried out during a visit of the first author
to the Ohio State University (supported in part by the MRI), and a
visit of the second author to the Universit\`a di Pisa. The work of
the first author was partly supported by the European Research
Council -- MEIF-CT-2005-010975 and MERG-CT-2007-046557. The work of
the second author was partly supported by NSF grants DMS-0606002, 
DMS-0906483, and by an Alfred P. Sloan research fellowship.

\vskip 10pt

\section{Background material on asymptotic cones.}

In this section, we provide some background on ultralimits and
asymptotic cones of metric spaces.  Let us start with some basic
reminders on ultrafilters.

\begin{Def}
A {\em non-principal ultrafilter} on the natural numbers $\mN$ is a
collection $\mU$ of subsets of $\mN$, satisfying the following four
axioms:
\begin{enumerate}
\item if $S\in \mU$, and $S^\prime \supset S$, then $S^\prime \in \mU$,
\item if $S \subset \mN$ is a finite subset, then $S\notin \mU$,
\item if $S,S^\prime \in \mU$, then $S\cap S^\prime \in \mU$,
\item given any finite partition $\mN=S_1\cup \ldots \cup S_k$ into pairwise
disjoint sets, there is a unique $S_i$ satisfying $S_i\in \mU$.
\end{enumerate}
\end{Def}

Zorn's Lemma guarantees the existence of non-principal
ultrafilters.  Now given a compact Hausdorff space $X$ and a map $f:
\mN\rightarrow X$, there is a unique point $f_\om \in X$ such
that every neighborhood $U$ of $f_\om$ satisfies $f^{-1}(U)\in
\mU$.  This point is called the $\om-$limit of the sequence $\{f(i)\}$;
we will occasionally write $\om\lim \{f(i)\} := f_\om$. 
In particular, if the target space $X$ is the compact space
$[0, \infty]$, we have that $f_\om$ is a well-defined real
number (or $\infty$).

\begin{Def}
Let $(X,d, *)$ be a pointed metric space, $X^\mN$ the collection of
$X$-valued sequences, and $\lambda:\mN\rightarrow (0,\infty) \subset
[0,\infty]$ a sequence of real numbers satisfying $\lambda_\om
=\infty$.  Given any pair of
points $\{x_i\}, \{y_i\}$ in $X^\mN$, we define the pseudo-distance
$d_\om(\{x_i\}, \{y_i\})$ between them to be $f_\om$, where
$f:\mN\rightarrow [0,\infty)$ is the function $f(k)=
d(x_k,y_k)/\lambda(k)$.  Observe that this pseudo-distance takes on
values in $[0,\infty]$.

Next, note that $X^\mN$ has a distinguished point,
corresponding to the constant sequence $\{*\}$.
Restricting to the subset of $X^\mN$
consisting of sequences $\{x_i\}$ satisfying $d_\om(\{x_i\},
\{*\})<\infty$, and identifying sequences whose $d_\om$ distance is
zero, one obtains a genuine pointed metric space $(X_\om, d_\om,
*_\om)$, which is called an {\em asymptotic cone} of the pointed metric
space $(X,d, *)$.
\end{Def}

We will usually denote an asymptotic cone by $Cone(X)$.  The reader
should keep in mind that the construction of $Cone(X)$ involves a
number of choices (basepoints, sequence $\lambda_i$, choice of
non-principal ultrafilters) and that different choices could give
different (non-homeomorphic) asymptotic cones (see the papers
\cite{TV}, \cite{KSTT}, \cite{OS}). However, in the special case where
$X=\mR ^k$, {\it all} asymptotic cones are isometric to $\mR ^k$ (i.e.
we have independence of all choices).

We will require the following facts concerning asymptotic cones of
non-positively curved spaces:
\begin{itemize}
\item if $(X,d)$ is a CAT(0)-space, then $Cone(X)$ is likewise a CAT(0)-space,
\item if $\phi: X\rightarrow Y$ is a $(C,K)$-quasi-isometric map, then $\phi$
induces a $C$-bi-Lipschitz map $\phi_\om: Cone(X)\rightarrow Cone(Y)$,
\item if $\gamma \subset X$ is a $k$-flat, then
  $\gamma_\omega:=Cone(\gamma)\subset
  Cone(X)$ is likewise a $k$-flat,
\item if $\{a_i\}, \{b_i\} \in Cone(X)$ are an arbitrary pair of points, then the ultralimit of the geodesic
segments $\overline{a_ib_i}$ gives a geodesic segment
$\overline{\{a_i\}\{b_i\}}$ joining $\{a_i\}$ to $\{b_i\}$.
\end{itemize}
Concerning the second point above, we remind the reader that a
$(C,K)$-quasi-isometric map $\phi: (X, d_X) \rightarrow (Y, d_Y)$ between
metric spaces is a (not necessarily continuous) map having the
property that:
$$\frac{1}{C} \cdot d_X(p,q) - K \leq d_Y(\phi(p), \phi(q))\leq C\cdot d_X(p,q) + K.$$
We also comment that, in the second point above, the asymptotic cones of $X$, $Y$,
have to be taken with the same scaling sequence and the same ultrafilters.

\begin{Lem}[Translations on asymptotic cone]\label{l_trans}
Let $X$ be a geodesic space, $\gamma \subset X$ a $k$-flat, and
$\gamma_\om \subset 
Cone(X)$ the corresponding $k$-flat in an asymptotic cone $Cone(X)$ of
$X$.  Assume that there exists a subgroup $G< Isom(X)$ with the
property that $G$ leaves $\gamma$  invariant, and acts cocompactly on
$\gamma$.  Then for any pair of points $p,q\in \gamma_\om$, 
there is an isometry $\Phi: Cone(X)\rightarrow Cone(X)$ satisfying $\Phi(p)=q$.
\end{Lem}

\begin{Prf}
Let $\{p_i\}, \{q_i\}\subset \gamma \subset X$ be sequences defining
the points $p,q$ respectively.  Since $G$ leaves $\gamma$ invariant,
and acts cocompactly on $\gamma$, there exists elements $g_i\in G$ 
with the property that for every index $i$, we have $d(g_i(p_i),
q_i)\leq R$. 

Now observe that the sequence $\{g_i\}$ of isometries of $X$
defines a self-map (defined componentwise) of the space
$X^{\mathbb N}$ of sequences of points in $X$.  Let us denote by
$g_\om$ this self-map, which we now proceed to show induces the
desired isometry on $Cone(X)$.  First note that it is immediate
that $g_\om$ preserves the pseudo-distance $d_\om$ on $X ^{\mathbb
N}$, and has the property that $d_\om(\{g_i(p_i)\}, \{q_i\})=0$. 
So to see that $g_\om$ descends to an isometry of
$Cone(X)$, all we have to establish is that for $\{x_i\}$ a
sequence satisfying $d_\om(\{x_i\}, *)<\infty$, the image sequence
also satisfies $d_\om(\{g_i(x_i)\}, *)<\infty$. But we have the
series of equivalences:
$$d_\om(\{x_i\}, *)<\infty \hskip 5pt \Longleftrightarrow \hskip 5pt d_\om(\{x_i\},\{p_i\})<\infty$$
$$\Longleftrightarrow \hskip 5pt d_\om(\{g_i(x_i)\},\{g_i(p_i)\})<\infty$$
$$\Longleftrightarrow \hskip 5pt d_\om(\{g_i(x_i)\},\{q_i\})<\infty$$
$$\Longleftrightarrow \hskip 5pt d_\om(\{g_i(x_i)\},*)<\infty$$
where the first and last equivalences come from applying the
triangle inequality in the pseudo-metric space $(X ^{\mathbb
N},d_\om)$, and the second and third equivalences follow from our
earlier comments. We conclude that the induced isometry $g_\om$ on
the pseudo-metric space $X^{\mathbb N}$ of sequences leaves
invariant the subset of sequences at finite distance from the
distinguished constant sequence, and hence descends to an isometry
of $Cone(X)$. Finally, it is immediate from the definition of the
isometry $g_\om$ that it will leave $\gamma_\om$ invariant, as each
$g_i$ leaves $\gamma$ invariant. This concludes the proof.
\flushright{$\square$}
\end{Prf}

Let us now specialize the previous Lemma to the case of geodesics (i.e. $k=1$).
Observe that any element $g\in Isom(X)$ as in the previous Lemma gives
rise to a $\mathbb Z$-action on $X$ leaving $\gamma$ invariant.
It is worth pointing out that the Lemma does {\bf not} state that the
$\mathbb Z$-action on $X$ induces an $\mathbb R$-action on $Cone(X)$. 
The issue is that for each $r\in \mathbb R$, there is indeed a
corresponding isometry of $Cone(X)$, but these will not in 
general vary continuously with respect to $r$ (as can already be
seen in the simple case where $X=\mathbb H^2$).

\section{Flattening sequences and Half-flats}
In this section, we will provide a proof of Theorem A. Our goal is to show how
certain sequences of maps from the disk to the asymptotic cone of a CAT$(0)$-space $M$
can be used to construct flats in $M$. We recall that $\omega$ denotes
the ultrafilter used to construct $X=Cone(M)$, that $\lambda_j$
denotes the sequence of scaling factors, and that $*$ denotes both the
base-point of $M$, and the basepoint of $X$ represented by the constant
sequence $\{*\}$.

We are given a $k$-flat $\gamma\subset M$ (possibly a geodesic), and we have the 
corresponding $k$-flat $\gamma_\omega \subset X$ in the asymptotic cone $X$. By abuse
of notation, we will use $\pi$ to denote both the nearest point projection $\pi: M \rightarrow \gamma$,
as well as the nearest point projection $\pi: X \rightarrow \gamma_\omega$. We can now make the:

\begin{Def}[Flattening sequences]
We say that $\gamma_\omega$ has a {\em flattening sequence} 
provided there exists a sequence of continuous maps $(f_r)_{r\in \mathbb N}$ from 
the $k$-disk $\mD^k$ to $X$ such that
\begin{enumerate}
\item The diameter of $f_r(\mD^k)$ is smaller than one.
\item The image of $\pi\circ f_r$ contains an open
  neighborhood of the base-point $*$.
\item The restriction of $\pi\circ f_r$ to $\partial \mD^k= S^{k-1}$ does not
  contain $*$, and represents a non-zero element in the homotopy group $\pi_{k-1}(\gamma_\omega\setminus \{*\}) \cong \mZ$.
\item There exists a constant $D>0$ with the property that
$$d(f_r(\mD^k),\gamma_\omega)= \inf_{x\in \mD^k}d(f_r(x), \gamma_\omega) \geq D\cdot r.$$
\end{enumerate}
The sequence of maps $(f_r)_{r\in \mN}$ will be called a {\em flattening sequence} for $\gamma_\omega$.
\end{Def}

We now assume that the $k$-flat $\gamma_\omega$ has a flattening sequence consisting
of maps $f_r: \mD^k \rightarrow X$. To establish Theorem A, we need to prove that $\gamma$
bounds a $(k+1)$-dimensional half-flat. In order to do this, we have to construct geodesic rays 
emanating from various points on $\gamma$. 

Each such ray will be constructed as a limit of a sequence
of longer and longer geodesic segments, originating from a fixed point on $\gamma$, and terminating
at a sequence of suitably chosen points in the space $M$. 
In order to select this suitable sequence of points in $M$, we start by noting that each of the maps 
$f_r: \mD ^k \rightarrow X$ in our flattening sequence can be obtained as an ultralimit of a sequence of
maps $f_{r,j}: \mD^k \rightarrow M$ (see for instance Kapovich \cite{Ka}). The desired collection of points 
will be carefully chosen  to lie on the image of some of the maps $f_{r,j}$. 
The precise selection process is contained in the following:

\vskip 10pt

\noindent {\bf Assertion:} Let us be given an arbitrary finite $(m+1)$-tuple of points $\{P^0, \ldots , P^m\} \subset \gamma$. 
Then for each $r\in \mN$, we can choose indices $j_r \in \mN$, and $(m+1)$-tuples
of points $\{x^0_{r,j_r}, \ldots , x^m_{r, j_r}\} \subset f_{r,j_r}\big(Int(\mD^k)\big) \subset M$ with the property that:
  \begin{enumerate}
  \item For each $r$, and $0\leq i \leq m$, we have that $\pi(x^i_{r,j_r})=P^i$.
  \item For any $i,i^\prime$ $$\frac{d(x^i_{r,j_r},x^{i^\prime}_{r,j_r})}{\lambda_{j_r}} < 2.$$
  \item For any $i$ $$\frac{d(P^i,x^i_{r,j_r})}{r \cdot \lambda_{j_r}}> D/2.$$
  \end{enumerate}

\vskip 10pt

We temporarily delay the proof of the {\bf Assertion}, and focus on explaining how our Theorem A can be 
deduced from this statement. We first recall some terminology: we
say that two half-rays $\eta_1$ and $\eta_2$ {\it bound a flat strip}, provided there exists
an isometric embedding of $\mathbb R^+ \times [0, a]$ (for some $a>0$) into $M$, with
the property that $\eta_1$ coincides with $\mathbb R^+\times \{0\}$ and $\eta_2$ coincides
with $\mathbb R^+\times \{a\}$. For a fixed ray $\eta$, we will denote by $Par(\eta)\subset M$ the union
of all geodesic rays in $M$ which, together with $\eta$, jointly bound a flat strip. Our first step
is to use the {\bf Assertion} to show:

\vskip 10pt

\noindent {\bf Claim 1:} Given any finite set of points $\{P^0, \ldots , P^m\} \subset \gamma$, there exist
a collection of geodesic rays $\{\eta ^0, \ldots , \eta ^m\}$, satisfying:
\begin{itemize}
\item for each $s$, we have that $\eta ^i$ originates at $X^i$, and satisfies $\pi(\eta^i) = P^i$.
\item each pair of geodesic rays $\eta^i, \eta ^{j}$ jointly bound a flat strip.
\end{itemize} 

\vskip 10pt

\begin{Prf}[Claim 1]
This can be seen as follows: first, apply the {\bf Assertion} to the finite set of points, obtaining a sequence
of $(m+1)$-tuples of points $\{x^0_{r,j_r}, \ldots , x^m_{r, j_r}\} \subset M$. Now for each $i$, consider the
sequence of geodesic segments $\eta ^i_r$, which joins the point $X^i$ to the point $x^i_{r,j_r}$. Since
$M$ is locally compact, we can extract a subsequence which simultaneously converges for all the $0\leq i\leq m$.
We define the limiting geodesic rays to be our $\eta ^i$. So to complete the proof of the Claim, we just need
to verify that these $\eta^i$ have the desired property. By construction, we know that each of the $\eta^i_r$ 
are geodesic segments originating at $P^i$, which immediately gives us the corresponding property for 
$\eta^i$. Likewise, we have that each of the geodesic segments $\eta^i_r$ project to the point $P^i$, 
which yields the same statement for $\eta ^i$. Note that this implies that, for any $t>0$, $\eta^i$ is the minimal
length path from $\eta^i(t)$ to $\gamma$. In particular, the angle (in the CAT(0) sense, see \cite{BrHa}) of 
$\eta^i$ with $\gamma$ must be $\geq \pi/2$.

So we now need to establish the second property: that $\eta ^i$ and $\eta ^{j}$ jointly bound a 
flat strip. Let us set $d= d(P^i, P^{j})$; we start by showing that for $t\geq 0$, we have 
$d\big(\eta^i(t), \eta^{j}(t)\big)=d$. Since the geodesic rays are limits of the geodesic segments, this is equivalent
to showing $\lim_{r\to \infty} d\big(\eta^i_r(t), \eta ^j_r(t)\big) =d$. From condition (1) in the 
{\bf Assertion}, and the fact that the projection map is distance non-increasing, we obtain the inequality:
$$\lim_{r\rightarrow \infty} d\big(\eta^i_r(t), \eta ^j_r(t)\big) \geq \lim_{r\rightarrow \infty} 
d\big(\eta^i_r(0), \eta ^j_r(0)\big) = d(P^i, P^j) =d$$
For the reverse inequality, we need to estimate from above the distance between $\eta^i_r(t)$ and 
$\eta ^j_r(t)$. To do this, we first truncate the longer of the two geodesic segments 
$\eta^i_r$ and $\eta^j_r$ to have length equal to the smaller one. We will denote by $\bar \eta^i_r$,
$\bar \eta^j_r$ the new pair of equal length geodesics, and let $L$ denote their common length.
Since we are trying to estimate from above the distance between the points $\bar \eta ^i_r(t)$ and 
$\bar \eta ^j_r(t)$, convexity of the distance function in CAT(0)-spaces yields:

$$d\big(\bar \eta ^i_r(t) , \bar \eta ^j_r(t)\big) \leq \Big(1-\frac{t}{L}\Big) \cdot d\big(\bar \eta ^i_r(0) , \bar \eta ^j_r(0)\big)
+ \frac{t}{L} d\big(\bar \eta ^i_r(L) , \bar \eta ^j_r(L)\big)$$

But from condition (3) in the {\bf Assertion}, we have that $L> r \cdot \lambda _{j_r} \cdot D/2$. Furthermore,
from condition (2) in the {\bf Assertion}, and an application of the triangle inequality, we obtain that 
$d\big(\bar \eta ^i_r(L) , \bar \eta ^j_r(L)\big) < 4 \lambda _{j_r} + d$. Substituting these expressions (along with
$d\big(\bar \eta ^i_r(0) , \bar \eta ^j_r(0)\big)=d$) into
our inequality, and simplifying, we get:
$$d\big(\bar \eta ^i_r(t) , \bar \eta ^j_r(t)\big) \leq d + 4t \cdot \frac{\lambda_{j_r}}{L} < d + 4t \cdot \frac{\lambda_{j_r}}
{r \lambda_{j_r} \cdot D/2} = d + \frac{8t}{D\cdot r}$$
where we recall that $d,t,D$ are constants.
Now taking the limit as $r\rightarrow \infty$, we obtain that $\lim_{r\to \infty} 
d\big(\bar \eta ^i_r(t) , \bar \eta ^j_r(t)\big) \leq d$, as desired.

This verifies that the two geodesic rays $\eta^i$, $\eta^j$ remain at a constant distance apart, in the sense that
$d\big(\eta^i(t), \eta^j(t)\big)$ is a constant function of $t$. Furthermore, from our earlier discussion, we
have that the geodesics segment joining $P^i= \eta^i(0)$ to $P^j=\eta^j(0)$ forms an angle 
$\geq \pi/2$ with each of the geodesic rays $\eta^i, \eta^j$, and hence both these angles must 
actually be $=\pi/2$. But in a CAT(0) space, this forces the 
geodesics $\eta^i$ and $\eta^j$ to bound a flat strip (see the proof of the flat strip theorem \cite[pg. 182]{BrHa}).
This concludes the proof of {\bf Claim 1}. \qed

\end{Prf}

So we now know that any finite set of points $\{P^0, \ldots ,P^m\} \subset \gamma$ are common endpoints of
geodesic rays that pairwise bound a flat strip. But ultimately, we want to show that {\it every} point on $\gamma$ 
is an endpoint of a parallel geodesic ray. Our next step is to establish:

\vskip 10pt

\noindent {\bf Claim 2:} Given any compact set $K\subset \gamma$, we can find an isometric embedding of
$K\times [0,\infty) \hookrightarrow M$ with the property that $K\times \{0\}$ maps to $K$.

\vskip 10pt

\begin{Prf}[Claim 2]
Recall that, for a geodesic $\eta$, the set $Par(\eta)$ is the union of all geodesic rays which,
together with $\eta$, bound a flat strip.
From the proof of the product region theorem in CAT(0)-spaces (see for example \cite{BrHa}), one has that
$Par(\eta)$ forms a convex subset of $M$, which splits as a metric  product $B\times [0,\infty)$. 
Here $B\subset M$ is a convex subset, and consists of the collection of all the basepoints of the parallel
geodesic rays. 

Now, since $K\subset \gamma$ is compact, we can find a finite set of points $\{P^0,\ldots ,P^m\} \subset \gamma$ 
whose convex
hull contains $K$. Applying {\bf Claim 1}, we have that there exist a corresponding collection of geodesic rays 
$\{\eta ^0, \ldots ,\eta ^m\}$, with each $\eta^i$ originating from $P^i$, and which pairwise bound a flat strip.
Considering the set $Par (\eta ^0)$, we have that $Par(\eta ^0)$ is isometric to $B\times [0,\infty)$, where 
$B$ is the collection of basepoints.
But we know that $\{P^0,\ldots, P^m\} \subset B$, so their convex hull is likewise contained in $B$, forcing
$K\subset B$. We conclude that there is an isometric copy of $K\times [0,\infty)$ embedded inside the convex
subset $Par(\eta) \subset M$. This completes the proof of {\bf Claim 2}. \qed
\end{Prf}

Finally, let us take a sequence of compact sets $K_i\subset \gamma$ exhausting the flat $\gamma$ (for instance,
take $K_i$ to be the radius $i$ metric ball centered at $*$). From {\bf Claim 2}, we have a corresponding 
sequence of isometrically embedded copies of $K_i \times [0,\infty) \hookrightarrow M$, where each $K_i\times \{0\}$
maps to the corresponding compact $K_i$. From local compactness, we can extract a convergent subsequence,
whose limit will be the desired half-flat bounding $\gamma$. So to complete the proof of Theorem A, we are left 
with verifying the {\bf Assertion}. 

\vskip 10pt

\begin{Prf}[Assertion] Let us recall the framework: we have a finite collection of points 
$\{P^0,\ldots , P^m\}\subset \gamma$, and for each $r\in \mN$, a collection of maps $f_{r,j}: \mD^k \rightarrow M$
having the property that $\omega \lim f_{r,j} = f_r : \mD^k \rightarrow X$. We want to find, for each index $r$, a
corresponding index $j_r$ and set of points $\{x^0_{r,j_r}, \ldots , x^m_{r,j_r}\} \subset f_{r,j_r}\big(Int(\mD^k)\big)$.
The chosen set of points should have the property that (1) for each $i$, $\pi(x^i_{r,j_r})=P^i$,
(2) for any $i,i^\prime$, ${d(x^i_{r,j_r},x^{i^\prime}_{r,j_r})} < 2{\lambda_{j_r}}$, and 
(3) for any $i$, ${d(P^i,x^i_{r,j_r})}> {r \cdot \lambda_{j_r}}\cdot D/2$.

Our approach is as follows: fixing $r$, we define three subsets of $\mN$ by setting
\begin{itemize}
\item $J_1$ to be the set of indices $j$ for which $\{P^0, \ldots, P^m\} \subset \pi \circ f_{r,j} \big( Int(\mD^k)\big)$ ,
\item $J_2$ to be the set of indices where $diam\big(f_{r,j}(\mD^k)\big)<2 \lambda_j$, and
\item $J_3$ to be the set of indices where $d\big(\gamma, f_{r,j}(\mD^k)\big) > {r \cdot \lambda_{j}}\cdot D/2$.

\end{itemize}
Now assuming we could show that each of these three sets are in $\omega$, property (3) of ultrafilters
(closure under finite intersections) implies that $J_1\cap J_2 \cap J_3 \in \omega$. 
Finally, every set in $\omega$ is infinite (property (2) of ultrafilters), and in particular non-empty, 
allowing us to find an index $j_r \in J_1\cap J_2\cap J_3$. For this index $j_r$, we can choose
arbitrary points $x^i_{r,j_r} \in f_{r,j_r}\big(Int(\mD^k)\big)$ satisfying $\pi(x^i_{r,j_r})=P^i$ (such
points exist since $j_r \in J_1$). And from $j_r\in J_2\cap J_3$, it immediately follows that
the tuple of points $\{x^0_{r,j_r}, \ldots , x^m_{r,j_r}\} \subset f_{r,j_r}\big(Int(\mD^k)\big)$ has the
desired properties.

\vskip 5pt

\noindent {\bf Step 1:} The set $J_2$ lies in the ultrafilter $\omega$.

\vskip 5pt

To see this, we recall that property (1) in the definition of a flattening sequence
requires $diam\big(f_r(\mD ^k)\big) \leq 1$. Since we know that $f_r = \omega \lim f_{r,j}$,
the definition of distances in the ultralimit tells us that the set of indices $j$ for which 
$diam\big(f_{r,j}(\mD^k)\big)/\lambda_j<2$ lies in the ultrafilter. This verifies Step 1.

\vskip 5pt

\noindent {\bf Step 2:} The set $J_3$ lies in the ultrafilter $\omega$.

\vskip 5pt

For this, we argue similarly. Recall that property (4) in the definition of a flattening sequence
requires the existence of a $D>0$ so that
$d(f_r(\mD^k),\gamma_\omega)= \inf_{x\in \mD^k}d(f_r(x), \gamma_\omega) \geq D\cdot r$. Since
$f_r$ is the ultralimit of the maps $f_{r,j}$, the definition of distances in the ultralimit tells us that 
the set of indices $j$ for which 
$d(f_{r,j}(\mD^k),\gamma)/\lambda_{j} \geq rD/2$ lies in the ultrafilter. This verifies Step 2.

\vskip 5pt

\noindent {\bf Step 3:} The set $J_1$ lies in the ultrafilter $\omega$.

\vskip 5pt

This last step is much more involved than the first two. Let us fix one of the points $P^i$, and 
consider the restriction $f_{r,j} | _{\partial D^k}: S^{k-1} \rightarrow M$, composed with the projection
$\pi: M \rightarrow \gamma$. We have three distinct possibilities:
\begin{enumerate}
\item $P^i$ lies in the image of $\pi \circ f_{r,j} | _{\partial D^k}$, or
\item $\pi\big(f_{r,j}(\partial D^k)\big) \subset \gamma -\{P^i\}$, and $\big[\pi \circ f_{r,j}|_{\partial D^k}\big]=0$ in $\pi_{k-1}(\gamma -\{P^i\})$, or
\item $\pi\big(f_{r,j}(\partial D^k)\big) \subset \gamma -\{P^i\}$, and $\big[\pi \circ f_{r,j}|_{\partial D^k}\big]\neq0$ in $\pi_{k-1}(\gamma -\{P^i\})$.
\end{enumerate}
This gives us a partition $\mN = I_1^i \cup I_2^i \cup I_3^i$ into three disjoint sets, according to which of these three
properties holds for the index $j$. From property (4) of ultrafilters, we have that exactly one
of these three sets must lie in $\omega$. If we could show that $I_3^i\in \omega$, then property (3) of ultrafilters
would force $I_3^0\cap \ldots \cap I_3^k \in \omega$. Since we have a containment 
$I_3^0\cap \ldots \cap I_3^k \subset J_1$, property (1) of ultrafilters would give us that $J_1 \in \omega$.

\vskip 5pt

So to conclude the proof of Step 3 (and hence, of the {\bf Assertion}), we are left with showing that $I_1^i\notin \omega$ and 
that $I_2^i \notin \omega$. We will argue both of these by contradiction. Supposing that $I_1^i \in \omega$, 
we would have that the set of indices for which 
$d(P^i , \pi \circ f_{r,j} | _{\partial D^k})/\lambda_j = 0$ is contained in $\omega$. So
the point represented by the constant sequence $\{ P^i \} = \{ *\} \in \gamma _\omega$ lies on
the set:
$$\omega \lim \big( \pi \circ f_{r,j} | _{\partial D^k} \big) = \pi \big ( \omega \lim f_{r,j} | _{\partial D^k}\big)
 = \pi \circ f_r | _{\partial D^k}$$
But this contradicts property (3) in the definition of flattening sequence. 

Similarly, to see that $I_2^i\notin \omega$, we again argue by contradiction. So let us assume that 
$I_2^i \in \omega$. Note that the indices in $I_2^i$ are those for which the map $\pi \circ f_{r,j}|_{\partial D^k}:
S^{k-1} \rightarrow \gamma -\{P^i\}$ is homotopically trivial. Let us define real numbers 
$$a_j=\inf _{x\in \partial D^k} d\big(P^i, \pi \circ f_{r,j}(x)\big),$$ 
$$b_j=\sup _{x\in \partial D^k} d\big(P^i, \pi \circ f_{r,j}(x)\big)$$
and consider the set $A_j:= \{ x\in \gamma \hskip 5pt | \hskip 5pt a_j \leq d(x, P_i) \leq b_j\}$. The inclusion 
$A_j \hookrightarrow \gamma - \{P^i\}$ is a homotopy equivalence, and the map $\pi \circ f_{r,j}|_{\partial D^k}$
has image lying inside $A_j$. Hence for indices $j\in I_2^i$, we have that 
$[\pi \circ f_{r,j}|_{\partial D^k}] = 0 \in \pi_{k-1}(A_j)$, and we can construct a map 
$F_j: \mD^k \rightarrow A_j \subset \gamma - \{P^i\}$
with the property that $F_j | _{\partial \mD^k} = \pi \circ f_{r,j}|_{\partial D^k}$. 

Now for each $\epsilon >0$,
we can further partition the set $I_2^i = I_2^i (\epsilon) \cup \bar I_2^i(\epsilon)$ by defining $I_2^i(\epsilon)$
to be the set of indices where the inequality $a_j < \epsilon \lambda_j$ holds, and $\bar I_2^i(\epsilon)$ 
the set of indices where $a_j \geq \epsilon \lambda_j$. From property (4) of ultrafilters,
we have that precisely one of the sets $I_2^i(\epsilon)$, $\bar I_2^i(\epsilon)$ is contained in $\omega$.
We now have two possibilities:
\begin{itemize}
\item either we have that $I_2^i(\epsilon) \in \omega$ for {\it every} $\epsilon >0$, or  
\item there exists some $\epsilon >0$, for which $\bar I_2^i(\epsilon) \in \omega$.
\end{itemize}
In the first case, we again obtain that the point represented by the constant sequence 
$\{P^i\}=\{*\} \in \gamma_\omega$ lies on the set $\omega \lim \big( \pi \circ f_{r,j} | _{\partial D^k} \big)
= \pi \circ f_r | _{\partial D^k}$, which contradicts property (3) in the definition of flattening sequence.

In the second case, we can try to take the ultralimit of the collection of maps 
$F_j : \mD^k \rightarrow \gamma$. The {\it upper bound} on the distance between $F_j$ and $P^i$ ensures
that the $F_j$ escape to infinity slowly enough for the ultralimit to be defined. More precisely,  
from the fact that the maps $\pi \circ f_{r,j}|_{\partial D^k}$ have as ultralimit $\pi \circ f_r| _{\partial \mD^k}$, 
we must have that $\omega \lim \{b_j /\lambda _j\} < \infty$, which
in turn implies that the ultralimit $F_\omega : \mD^k \rightarrow \gamma_\omega$ exists. On the other
hand, the {\it lower bound} on the distance between $F_j$ and $P^i$ ensures that the ultralimit $F_\omega$
does not pass through the constant sequence $\{P^i\}$. More precisely, for the $\epsilon>0$  satisfying 
$\bar I_2^i(\epsilon)\in \omega$, and for any index set $j\in \bar I_2^i(\epsilon)$, we have that
$$d\big( F_j(\mD^k) , P^i\big)/\lambda _j \geq d( A_j , P^i)/\lambda _j = a_j/\lambda _j \geq \epsilon.$$
Since this holds for a set of indices in the ultrafilter, we immediately deduce the corresponding property
for the ultralimit: $d\big( F_\omega (\mD^k) , \{P^i\} \big) \geq \epsilon$. In particular, $F_\omega:\mD^k
\rightarrow \gamma_\omega$
has image lying in the complement of the point corresponding to the constant sequence $\{P^i\}=\{*\}$, 
and restricts to the map
$\pi \circ f_r |_{\partial \mD^k}$ on the boundary $S^{k-1}= \partial \mD^k$. This tells us that 
$[\pi \circ f_r |_{\partial \mD^k}] = 0 \in \pi_{k-1}(\gamma_\omega \setminus \{*\})$, which contradicts property (3)
in the definition of a flattening sequence. This concludes the verification of Step 3, and hence completes
the proof of the {\bf Assertion}.
\qed
\end{Prf}

\vskip 5pt

\noindent Having established the {\bf Assertion}, we have now concluded the proof of Theorem A.
From the definition of a flattening sequence, it is obvious that these exist whenever $\gamma_\omega$
bounds a half-flat in $X$. As a result, we see that Theorem B follows immediately from Theorem A. 

\vskip 5pt

Finally, let us conclude this section by providing a family of cautionary examples. These will be 
locally compact CAT(0)-spaces $X_k$, each containing a geodesic $\gamma$, with the property that for
a suitable choice of scales, $\gamma_\omega$ is contained inside a $k$-dimensional flat, but the
individual $X_k$ do not contain any flats of dimension $>1$. In particular, these examples show that
the analogue of Theorem B with ``half-flats'' replaced by ``flats'' is {\it false}.

\vskip 10pt

\noindent {\bf Example:} Let us fix a $k\geq 2$, and for $n\in \mN$, define the spaces 
$C_n := [-n^3, n^3]^{k} \subset \mR^k$. 
Each $C_n$ is isometric to the standard $k$-dimensional cube with side lengths $2n^3$; we let 
$l_n \subset C_n$ 
be the geodesic segment of length $2n^3$ joining the two points $(\pm n^3, 0, \ldots, 0)$ 
inside $C_n$. Now consider 
the closed upper half space $\mR \times \mR_{\geq 0}:=\{(x,y) \hskip 5pt | \hskip 5pt y\geq 0\}$ with the
standard flat metric, and for $n\in \mN$, let us denote by 
$\bar l_n$ the segment of length $2n^3$ joining the pair of points $(\pm n^3, n)$ inside 
$\mR\times \mR_{\geq 0}$. We now form the space $X_k$
by gluing together all the $C_n$ to $\mR\times \mR_{\geq 0}$. More precisely, we isometrically 
identify each $l_n$ with the corresponding $\bar l_n$.

Observe that this space $X_k$, with the natural induced metric, is a locally compact CAT(0)-space. 
Furthermore, it is clear that $X_k$ {\it does not contain any flats of dimension $>1$}. Now consider the ultralimit $X$
obtained by fixing the origin as the sequence of basepoints, and setting $\lambda (i)=i^2$ to be the
sequence of scales. Let us consider the geodesic $\gamma\subset X_k$ given by the $x$-axis in the 
$\mR\times \mR_{\geq 0}$ portion of $M$. We claim that the corresponding geodesic 
$\gamma_\omega \subset \omega\lim X_k$ 
is contained inside a $k$-flat in $\omega \lim X_k$. Indeed, this follows readily from the following two observations
\begin{itemize}
\item since the distance from $\bar l_n$ to the $\gamma$ grows linearly, while the scaling factor $\lambda$
grows quadratically, every point $P\in \gamma_\omega$ can be represented by a sequence $\{p_n\}$ with 
the additional property that $p_n \in C_n$, and
\item since the size of the cubes $C_n$ grows cubically, while the scaling factor $\lambda$ grows quadratically,
the subset $C_\omega \subset \omega \lim X_k$ consisting of all points having a representative sequence of the form 
$\{c_i\}$ (with each $c_i\in C_i$) is isometric to the standard $\mR^k$.
\end{itemize}
This concludes our family of locally CAT(0) examples.

\vskip 10pt

\section{From bi-Lipschitz half-ultraflats to half-flats}

In this section we prove Theorem C, allowing us to deduce the presence of half-flats in $M$
from the presence of bi-Lipschitz half-flats in the ultralimit $X$ along with a mild periodicity
condition. 

\vskip5pt

The context is the following: we have a locally compact 
CAT$(0)$-space $M$ (for instance, the universal cover of a
non-positively curved Riemannian manifold) and an
asymptotic cone $X$ of $M$. 
We have a $k$-flat $\gamma$ in $M$, and its limit $\gamma_\omega$ in
$X$. Moreover, we are supposing  that there exists $G<Isom(M)$ that
acts co-compactly on $\gamma$.  We are assuming that there is a bi-Lipschitz embedding 
$\phi:\mathbb R^k\times \mathbb R_{\geq 0} \to X$, 
whose restriction to $\mathbb R^k\times \{0\}$ maps onto 
$\gamma_\omega$, and we want to show that $\gamma$ 
bounds a $(k+1)$-dimensional half-flat in $M$. In view of our Theorem A, it is sufficient to 
find a flattening sequence for $\gamma_\omega$.

\vskip5pt

Let $C$ be the bi-Lipschitz
constant of $\phi$, and let us make the simplifying assumption that the restriction of $\phi$ 
to the set $\mathbb R^k \times \{0\}$ is actually an {\it isometry} onto $\gamma_\omega$; we
will use this isometry to identify $\gamma_\omega$ with $\mR ^k\times \{0\}$. 
 For $r\in \mR$, let us denote by $L_r = \mR^k\times\{r\} \subset \mR^k \times \mR_{\geq 0}$
the horizontal flat at height $r$. We will use $\rho$ to denote the obvious projection map
$\rho: L_r \rightarrow L_0$. To make our various expressions
more readable, we use $d$ to denote
distance in $X$ (as opposed to $d_\om$), and the norm
notation to denote distance inside $\mR ^k\times\mR_{\geq 0}$.  

We now define, for each $r\in [0,\infty)$ a map $$\psi_r: L_r
\rightarrow L_0$$ as follows: given $p\in L_r$, we have $\phi(p)\in
X$.  Since $\gamma_\om \subset X$ is a flat
inside the CAT(0) space $X$, there is a well defined,
distance non-increasing, projection map $\pi: X\rightarrow
\gamma_\om$, which sends any given point in $X$ to the
(unique) closest point on $\gamma_\om$. Hence, given $p\in L_r$, we
have the composite map $\pi\circ \phi: L_r\rightarrow \gamma_\om.$
But recall that, by hypothesis, $\phi$ maps $L_0$ homeomorphically to
$\gamma_\om$. We can now set $\psi_r: L_r\rightarrow L_0$ to be the
composite map $$\psi_r=\phi^{-1}\circ \pi\circ \phi$$ We now show that
$\psi_r$ is at finite distance from the projection map $\rho: L_r\rightarrow L_0$.

We first observe
that for arbitrary $x \in L_r$, we have that the distance from $x$ to $L_0$
is exactly $r$, and hence
from the bi-Lipschitz estimate, we have
$$d(\phi(x), \gamma_\om)=d(\phi(x), \phi(L_0))\leq Cr$$
Since $\pi$ is the nearest point projection onto $\gamma_\om$, this implies that 
$d\big(\phi(x), (\pi\circ \phi)(x)\big)\leq Cr$.
Since $(\pi\circ \phi)(x) = \phi ( \psi_r(x))$, we can again use the
bi-Lipschitz estimate to conclude that:
$$Cr \geq d\big(\phi(x), (\pi\circ \phi)(x)\big)=d(\phi(x), \phi(\psi_r(x))\geq \frac{1}{C} 
\cdot ||x - \psi_r(x)||$$
Which gives us the estimate $ ||x - \psi_r(x)||\leq C^2r$. This implies that $\psi$ is at bounded
distance from the projection map $\rho$. Since the latter is a homeomorphism onto $L_0\cong \mR^k$, 
it follows that $\psi_r$ is surjective. 

Since $\psi_r$ is a surjective, Lipschitz map, its differential exists almost everywhere and
it is almost everywhere non-degenerate (see for example \cite[Chapter 3]{We}).
It follows that we can find a $k$-disk $D_r$ in $L_r$, of diameter smaller
than $1/C$, and a point $p_r$ in $\psi_r(D_r)$ such that
$\psi_r(\partial D_r)$ is homotopically non 
trivial in $L_0\setminus p_r$. By Lemma~\ref{l_trans}, we can find an
isometry $g_r$ of $X$, that leaves $\gamma_\omega$ invariant, and satisfies $g_r(\phi(p_r))=*$.

Now we define the maps $f_r: \mD^k \cong D_r\to X$ via the composition
$$f_r=g_r\circ\phi|_{D_r}.$$
and observe that we have $\pi\circ f_r=g_r\circ\phi\circ\psi_r$.
Moreover, since the diameter of $D_k$ is smaller than $1/C$, the
diameter of $f_r(D_k)$ is smaller than one. Finally, it is clear that
our choices for $D_k$ implies that $f_r$ satisfies all the conditions
for being a flattening sequence for $\gamma$. Invoking Theorem A
completes the proof of Theorem C in the special case where $\phi |_{\mR^k\times \{0\}}$
maps $\mR^k\times \{0\}$ isometrically onto $\gamma_\omega$.

Finally, it is easy to see that the general case can be deduced from this special case.
Indeed, given our arbitrary bi-Lipschitz embedding $\phi: \mR^k \times \mR _{\geq 0} \rightarrow
X$, let $\phi_0:\mR^k\times \{0\} \rightarrow \gamma_\omega$ be the restriction to $\mR^k\times \{0\}$.
We are assuming $\phi_0$ is onto, hence is a homeomorphism. Now consider the reparametrized 
map $\hat \phi: \mR^k \times \mR _{\geq 0}\rightarrow X$ obtained by taking the composition:
$$\xymatrix{
\mR^k\times \mR _{\geq 0} \ar[rr]^{\phi_0^{-1}\times Id} & & \mR^k\times \mR_{\geq 0} \ar[rr]^{\phi} & & X \\
}$$
As a composition of bi-Lipschitz maps, $\hat \phi$ is bi-Lipschitz, and by construction, we see that 
$\hat \phi$ maps $\mR^k\times \{0\}$ isometrically onto $\gamma_\omega$. This reduces the general
case to the special case we discussed above, completing the proof of Theorem C.
\qed

\section{Some applications}

Finally, let us discuss some consequences of our main results. As a first application, we obtain
some constraints on the behavior of a quasi-isometry between locally compact CAT(0)-spaces.

\begin{Cor}[Constraints on quasi-isometries]
Let $\tM_1,\tM_2$ be two locally compact CAT(0)-spaces, and assume that 
$\phi: \tM_1\rightarrow \tM_2$ is a quasi-isometry.  Let
$\gamma \subset \tM_1$ be a $k$-flat, $\gamma_\om\subset X_1:=Cone(\tM_1)$ the corresponding
$k$-flat in the asymptotic cone, and assume that there exists a bi-Lipschitz $(k+1)$-dimensional half-flat 
$F\subset X_1$ bounding the $k$-flat $\gamma_\om \subset X_1$.  Then we have the following dichotomy:
\begin{enumerate}
\item Non-periodicity: every $k$-flat $\eta$ at bounded distance from $\phi(\gamma)$ has the property that
$\eta/ Stab_{G}(\eta)$ is non-compact, where $G=Isom(\tM_2)$, or 
\item Bounding: every $k$-flat $\eta$ at bounded distance from $\phi(\gamma)$ bounds a $(k+1)$-dimensional 
half-flat.
\end{enumerate}
\end{Cor}

\begin{Prf}
This follows immediately from our Theorem C.  Assume that the first possibility
does not occur, i.e. there exists a $k$-flat $\eta$ at bounded distance from $\phi(\gamma)$
with the property that $Stab_G(\eta) \subset G= Isom(\tM_2)$ acts cocompactly on $\eta$.  
Now recall that the quasi-isometry $\phi:
\tM_1\rightarrow \tM_2$ induces a bi-Lipschitz homeomorphism $\phi_\om: Cone(\tM_1)
\rightarrow Cone(\tM_2)$.  Since $\eta \subset \tM_2$ was a $k$-flat at finite distance 
from $\phi(\gamma)$, we have the containment:
$$\phi_\om(\gamma_\om) \subseteq \eta _\om \subset Cone(\tM_2).$$ 
Since $\phi_\om(\gamma_\om)$ is a bi-Lipschitz copy of $\mR ^k$ inside the $k$-flat
$\eta_\om$, we conclude that $\phi(\gamma_\om)= \eta_\om$.  
But recall that we assumed that $\gamma_\om$ was contained inside a 
bi-Lipschitz flat $\gamma_\om \subset F \subset Cone(\tM_1)$, and hence we see that
$\eta_\om\subset \phi_\om(F)$ is likewise contained inside a bi-Lipschitz flat.
Since the hypotheses of Theorem C are satisfied, we conclude that $\eta$ must bound
a $(k+1)$-dimensional half-flat, concluding the proof of Corollary 5.1. \flushright{$\square$}

\end{Prf}

\vskip 5pt

The statement of our first corollary might seem somewhat complicated.  We 
now proceed to isolate a special case of most interest: 

\begin{Cor}[Constraints on perturbations of metrics]
Assume that $(M,g_0)$ is a closed Riemannian manifold of non-positive sectional curvature,
and assume that $N^k\hookrightarrow M$ is an isometrically embedded compact 
flat $k$-manifold with image $\gamma_0$.  
Let $ \tilde \gamma_0 \subset \tM$ be the $k$-flat obtained by taking a connected 
lift of $\gamma_0$, and assume that $\tilde \gamma_0$ bounds a 
$(k+1)$-dimensional half-flat $F_0$.  

Then if $(M,g)$ is any other Riemannian metric
on $M$ with non-positive sectional curvature, and $\gamma \subset M$ is an isometrically embedded
flat $k$-manifold (in
the $g$-metric) freely homotopic to $\gamma_0$, then the lift $\tilde \gamma \subset (\tM, \tilde g)$ must also
bound a $(k+1)$-dimensional half-flat $F$.
\end{Cor}

We can think of Corollary 5.2 as a ``non-periodic'' version
of the Flat Torus theorem.  Indeed, in the case where $F$ is $\pi_1(M)$-periodic, 
the Flat Torus theorem applied to $(M,g)$ implies that $\tilde \gamma$ is likewise contained
in a periodic flat.

\begin{Prf}
Since $M$ is compact, the identity map provides a quasi-isometry $\phi: (\tM, \tilde g_0)
\rightarrow (\tM, \tilde g)$.  The half-flat $F_0$ containing $\tilde \gamma_0$ gives rise to a flat $(F_0)_\om
\subset Cone(\tM, \tilde g_0)$ containing $(\tilde \gamma_0)_\om$.  In particular, we can apply
the previous Corollary 5.1.

Next note that, since $\gamma_0, \gamma$ are freely homotopic
to each other, there is a lift $\tilde \gamma$ of $\gamma$ which is at finite distance (in the 
$g$-metric) from the given $\tilde \gamma_0\subset (\tM, \tilde g)$.  Indeed, taking the free 
homotopy 
$H: N^k\times [0,1]\rightarrow M$ between $H_0=\gamma_0$ and $H_1=\gamma$, we 
can then take a lift $\tilde H: \mathbb R^k \times [0,1] \rightarrow \tM$ satisfying the initial 
condition $\tilde H_0 =\tilde \gamma_0$ (the given lift of $\gamma_0$).  
The time one map $\tilde H_1: \mathbb R^k \rightarrow \tM$ will 
be a lift of $H_1=\gamma$, hence a $k$-flat in $(\tM, \tilde g)$.  Furthermore, the distance 
(in the $g$-metric) between $\tilde \gamma_0$ and $\tilde \gamma$ will clearly be bounded
above by the supremum of the $g$-lengths of the (compact) family of maps 
$H_p:[0,1] \rightarrow (M,g)$, $p\in N^k$, defined by $H_p(t) = H(p,t)$.

Now observe that by construction, the $\tilde \gamma \subset (\tM, \tilde g)$ from the previous
paragraph has 
$Stab_G(\tilde \gamma)$ acting cocompactly on $\tilde \gamma$, where $G=Isom(\tM, \tilde g)$.
Hence the first possibility in the conclusion of Corollary 5.1 cannot occur, and we conclude
that $\tilde \gamma$ must bound a $(k+1)$-dimensional half-flat $F$, as desired.  This concludes the proof of
Corollary 5.2.  \flushright{$\square$}

\end{Prf}

\vskip 5pt

Next we recall some terminology from differential geometry: for $M$ a complete, simply connected, Riemannian
manifold of non-positive sectional curvature, the {\it rank} of a geodesic $\gamma \subset X$ is the dimension 
$rk(\gamma)$ 
of the vector space of parallel Jacobi fields along $\gamma$. Note that the unit tangent vector field is always
parallel, hence the rank of a geodesic is always $\geq 1$; a geodesic is said to have {\it higher rank} provided
$rk(\gamma)\geq 2$. A geodesic $\gamma$ that bounds a $2$-dimensional half-flat automatically has 
$rk(\gamma) \geq 2$, as the unit normal vector field within the half-flat will be a parallel Jacobi field along $\gamma$.
Finally, the {\it manifold} $M$ is said to have higher rank provided {\it every} geodesic $\gamma \subset M$
satisfies $rk(\gamma) \geq 2$. The celebrated rank-rigidity theorem, established independently by Ballmann \cite{Ba2}
and Burns-Spatzier \cite{BuSp}, states that if $M$ has higher rank then either (1) $M$ is isometric to an irreducible,
higher-rank symmetric space of non-compact type, or (2) $M$ is reducible, and splits isometrically as a product 
$M_1\times M_2$ of lower dimensional
manifolds of non-positive sectional curvature. Our next two applications will exploit the combination of our main 
results with the rank-rigidity theorem to deduce some information concerning manifolds of non-positive sectional
curvature.

\vskip 5pt

Now recall that the classic de Rham theorem \cite{dR} states that any simply connnected, 
complete Riemannian manifold admits a decomposition as a metric product $\tM=\mR^k \times 
M_1\times \ldots \times M_k$, where $\mR^k$ is a Euclidean
space equipped with the standard metric, and each $M_i$ is metrically irreducible (and not
$\mR$ or a point).  Furthermore, this decomposition is unique up to permutation of the factors.  
This result was recently generalized by Foertsch-Lytchak \cite{FoLy} to cover 
finite dimensional geodesic metric spaces (such as ultralimits of Riemannian
manifolds).  Our next corollary shows that, in the
presence of non-positive Riemannian curvature, there is
a strong relationship between splittings of $\tM$ and splittings of $Cone(\tM)$.

\begin{Cor}[Asymptotic cones detect splittings]
Let $M$ be a closed Riemannian manifold of non-positive curvature, $\tM$ the universal
cover of $M$ with induced Riemannian metric, and $X=Cone(\tM)$ an arbitrary asymptotic
cone of $\tM$.  If $\tM= \mR^k\times M_1\times \ldots \times M_n$ is the de Rham splitting of $\tM$ into irreducible factors, and $X= \mR ^l\times X_1\times \ldots \times X_m$ is the Foertsch-Lytchak splitting of $X$ into irreducible factors, then $k=l$, $n=m$, and up to a relabeling of the index set, we have that each $X_i=Cone(M_i)$.
\end{Cor}

\begin{Prf}
Let us first assume that $\tM$ is irreducible (i.e. k=0, n=1), and show that 
$X=Cone(\tM)$ is also irreducible (i.e. l=0, m=1).  By way of contradiction, let us assume that $X$
splits as a metric product, and observe that this clearly implies that every geodesic $\gamma 
\subset X$ is contained inside a flat. In particular, from our Theorem A, we see that
every geodesic inside $\tM$ must bound a $2$-dimensional half-flat, and hence must have higher rank. 
Applying the Ballmann, Burns-Spatzier
rank rigidity result, and recalling that $\tM$ was irreducible, we conclude that $\tM$
is in fact an irreducible higher rank symmetric space. But now Kleiner-Leeb have shown
that for such spaces, the asymptotic cone is irreducible (see \cite[Section 6]{KlL}), giving
us the desired contradiction.

Let us now proceed to the general case: from the metric splitting of $\tM$, we get a corresponding
metric splitting $Cone(\tM)=\mR^k \times Y_1\times \ldots \times Y_n$, where each $Y_i=Cone(M_i)$.  Since each $M_i$ is irreducible, the previous paragraph tells us that each $Y_i$ is likewise
irreducible.  So we now have two product decompositions of $Cone(\tM)$ into irreducible factors.  
So assuming that each $Y_i$ is distinct from a point and is not isometric to $\mR$, we 
could appeal to the uniqueness portion of Foertsch-Lytchak \cite[Theorem 1.1]{FoLy} to conclude
that, up to relabeling of the index set, each $X_i=Y_i=Cone(M_i)$, and that the Euclidean
factors have to have the same dimension $k=l$. 

To conclude the proof of our Corollary, we establish that if $M$ is a simply connected, complete,
Riemannian manifold of non-positive sectional curvature, and $\dim (M)\geq 2$, then $Cone(M)$
is distinct from a point or $\mR$.  First, recall that taking an arbitrary geodesic $\gamma 
\subset M$ (which we may assume passes through the basepoint $*\in M$), we get 
a corresponding geodesic $\gamma_\om \subset Cone(M)$, i.e. an isometric
embedding of $\mR$ into $Cone(M)$.  In particular, we see that $\dim (Cone(M))>0$.  To see 
that $Cone(M)$ is distinct from $\mR$, it is enough to establish the existence of three points 
$p_1,p_2,p_3\in Cone(M)$ such that for each index $j$ we have:
\begin{equation}
d_\om(p_j, p_{j+2}) \neq d_\om(p_j, p_{j+1})+ d_\om(p_{j+1}, p_{j+2})
\end{equation}
But this is easy to do: take $p_1,p_2$ to be the two distinct points on the geodesic 
$\gamma_\om$ at distance one from the basepoint $*\in Cone(M)$, so that $d_\om (p_1,p_2)=2$.  Observe that one can
represent the points $p_1,p_2$ via the sequences of points $\{x_i\}$, $\{y_i\}$ along $\gamma$ 
having the property that $*\in \overline{x_iy_i}$, and $d(x_i, *)= \lambda _i = d(*, y_i)$, where
$\lambda _i$ is the sequence of scales used in forming the asymptotic cone $Cone(M)$.  Now
since $\dim (M)\geq 2$, we can find another geodesic $\eta$ through the basepoint 
$*\in M$, with the property that $\eta \perp \gamma$.  Taking the sequence $\{z_i\}$ to lie on
$\eta$, and satisfy $d(z_i, *)=\lambda _i$, it is easy to see that this sequence defines a third point
$p_3\in Cone(M)$ satisfying $d_\om(p_3, *)=1$.  From the triangle inequality, we immediately
have that $d_\om(p_1,p_3)\leq 2$ and $d_\om(p_2,p_3)\leq 2$.  On the other hand, since the
Riemannian manifold $M$ has non-positive sectional curvature, we can apply Toponogov's theorem
to each of the triangles $\{*, x_i, z_i\}$: since we have a right angle at the vertex $*$, and we have
$d(*,x_i)=d(*,z_i)=\lambda_i$, Toponogov tells us that
$d(x_i,z_i)\geq \sqrt{2}\cdot \lambda_i$.  Passing to the asymptotic cone, this gives the lower
bound $d(p_1,p_3)\geq \sqrt{2}$, and an identical argument gives the estimate $d(p_2,p_3)\geq
\sqrt{2}$.  It is now easy to verify that the three points $p_1,p_2,p_3$ satisfy (11), and hence
$Cone(M) \neq \mR$, as desired.  This concludes the proof of Corollary 5.3.

\flushright{$\square$}
\end{Prf}

Before stating
our next result, we recall that the celebrated rank rigidity theorem of Ballmann, Burns-Spatzier 
was motivated by Gromov's
well-known rigidity theorem, the proof of which appears in the book
\cite{BGS}.  Our next corollary shows how in fact Gromov's rigidity theorem can 
now be directly deduced from the rank rigidity theorem.  This is our last:

\begin{Cor}[Gromov's higher rank rigidity \cite{BGS}]
Let $M^*$ be a compact locally symmetric space of $\mR$-rank $\geq 2$, with universal
cover $\tM ^*$ irreducible,
and let $M$ be a compact Riemannian manifold with sectional curvature $K\leq 0$.
If $\pi_1(M) \cong \pi_1(M^*)$, then $M$ is isometric to $M^*$, provided $Vol(M)=Vol(M^*)$.
\end{Cor}

\begin{Prf}
Since both $M$ and $M^*$ are compact with isomorphic fundamental groups,
the Milnor-\v Svarc theorem gives us quasi-isometries:
$$\tM^* \simeq \pi_1(M^*) \simeq \pi_1(M) \simeq \tM$$
which induce a bi-Lipschitz homeomorphism $\phi: Cone(\tM^*)\rightarrow
Cone(\tM)$.  Now in order to apply the rank rigidity theorem, we
need to establish that every geodesic in $\tM$ has rank $\geq 2$.  

We first observe that the proof of Corollary 5.2 extends almost 
verbatim to the present setting.  Indeed, in Corollary 5.2, we used
the identity map to induce a bi-Lipschitz homeomorphism between the 
asymptotic cones, and then appealed to Corollary 5.1.  The sole difference
in our present context is that, rather than using the identity map, we use the 
quasi-isometry between $\tM$ and $\tM^*$ induced by the isomorphism 
$\pi_1(M)\cong \pi_1(M^*)$.  This in turn induces a bi-Lipschitz homeomorphism
between asymptotic cones (see Section 2).  The reader can easily verify that the rest
of the argument in Corollary 5.2 extends to our present setting, establishing 
that {\it every lift to $\tM$ of a periodic geodesic} in $M$ has rank $\geq 2$.

So we now move to the general case, and explain why {\it every} geodesic in 
$\tM$ has higher rank.  To see this, assume by way of contradiction that there 
is a geodesic 
$\eta \subset \tM$ with $rk(\eta)=1$. Note that the geodesic $\eta$ cannot
bound a half-plane.  But once we have the existence of such an $\eta$,
we can appeal to results of Ballmann \cite[Theorem 2.13]{Ba1},
which imply that $\eta$ can be approximated (uniformly on compacts) 
by lifts of periodic 
geodesics in $M$; let $\{\tilde \gamma_i\}\rightarrow \eta$ be such an
approximating sequence.  Since each $\tilde \gamma _i$ has $rk(\tilde 
\gamma_i)\geq 2$, it supports a parallel Jacobi field $J_i$, which
can be taken to satisfy $||J_i|| \equiv 1$ and $\langle J_i, \tilde \gamma_i^\prime
\rangle \equiv 0$. Now we see that:
\begin{itemize}
\item the limiting vector field $J$ defined along $\eta$ exists, due to
the control on $||J_i||$,
\item the vector field $J$ along $\eta$ is a parallel Jacobi 
field, since both the ``parallel'' and ``Jacobi'' condition can be encoded
by differential equations with smooth coefficients, solutions to which 
will vary continuously with respect to initial conditions, and
\item $J$ will have unit length and will be orthogonal to $\eta^\prime$,
from the corresponding condition on the $J_i$.
\end{itemize}
But this contradicts our assumption that $rk(\eta)=1$.  
So we conclude that every geodesic $\eta \subset \tM$ must satisy $rk(\eta)\geq 2$, 
as desired.

From the rank rigidity theorem, we now obtain
that $\tM$ either splits as a product, or is isometric to an irreducible higher rank
symmetric space.  Since the asymptotic cone of the irreducible 
higher rank symmetric space is
topologically irreducible (see \cite[Section 6]{KlL}), and $Cone(\tM)$ is homeomorphic
to $Cone(\tM^*)$, we have that $\tM$ cannot split
as a product.  Finally, we see that $\pi_1(M) \cong \pi_1(M^*)$ acts cocompactly,
isometrically on two irreducible higher rank symmetric spaces $\tM$ and $\tM^*$.  
By Mostow rigidity \cite{Mo}, we have that the quotient spaces are, after suitably rescaling, 
isometric.  This completes our proof of Gromov's higher rank rigidity theorem. \flushright{$\square$}

\end{Prf}

\vskip 5pt

Finally, let us conclude our paper with a few comments on this last corollary.

\vskip 10pt

\noindent {\bf Remarks: (1)} The actual statement of Gromov's theorem in \cite[pg. (i)]{BGS} 
does not assume $\tM^*$ to be irreducible, but rather $M^*$ to be irreducible (i.e. there
is no finite cover of $M^*$ that splits isometrically as a product).  This leaves the possibility
that the universal cover $\tM^*$ splits isometrically as a product, but no finite cover of $M^*$
splits isometrically as a product.  However, in this specific case, the desired result was already 
proved by Eberlein (see \cite{Eb}).  And in fact, in the original proof of Gromov's rigidity theorem,
the very first step (see \cite[pg. 154]{BGS}) consists of appealing to Eberlein's result
to reduce to the case where $\tM^*$ is irreducible.

\noindent {\bf (2)} In the course of writing this paper, the authors learnt of the 
existence of another proof of Gromov's rigidity result, which bears some similarity
to our reasoning. As the reader has surmised from
the proof of Corollary 5.4, the key is to somehow show that $M$ also has to have
higher rank. But a sophisticated result of Ballmann-Eberlein \cite{BaEb} establishes 
that the rank of a non-positively curved Riemannian manifold $M$ can in fact 
be detected directly from algebraic properties of $\pi_1(M)$, and hence the property
of having ``higher rank'' is in fact algebraic (see also the recent preprint of 
Bestvina-Fujiwara \cite{BeFu}).  The main advantage
of our approach is that one can deduce Gromov's rigidity result directly from
rank rigidity, and indeed, that one can geometrically ``see'' that the property of having higher
rank is preserved.

\noindent {\bf (3)} We point out that various 
other mathematicians have obtained results extending
Gromov's theorem (and which do not seem tractable using our
methods).  A variation considered by Davis-Okun-Zheng (\cite{DOZ}, 
requires $\tM^*$ to be {\it reducible} and
$M^*$ to be irreducible (the same hypothesis as in Eberlein's rigidity result).
However, Davis-Okun-Zheng allow the metric on $M$ to be locally
CAT(0) (rather than Riemannian non-positively curved), and are still
able to conclude that $M$ is isometric (after rescaling) to $M^*$. The 
optimal result in this direction is due to Leeb \cite{L}, giving a characterization
of certain higher rank symmetric spaces and Euclidean buildings within
the broadest possible class of metric spaces, the Hadamard spaces (complete
geodesic spaces for which the distance function between pairs of geodesics
is always convex).  It is worth mentioning that Leeb's result relies heavily
on the viewpoint developed in the Kleiner-Leeb paper \cite{KlL}. 

\noindent {\bf (4)} We note that our method of proof 
can also be used to establish a {\it non-compact, finite volume} analogue 
of the previous corollary. 
Three of the key ingredients going into our proof were (i) Ballmann's 
result on the density of periodic geodesics in the tangent bundle, (ii) 
Ballmann-Burns-Spatzier's rank rigidity theorem, and (iii) Mostow's
strong rigidity theorem.  A finite volume version
of (i) was obtained by Croke-Eberlein-Kleiner (see
\cite[Appendix]{CEK}), under the assumption that 
the fundamental group is finitely generated. A finite volume version
of (ii) was obtained by Eberlein-Heber (see \cite{EbH}). The finite
volume versions of Mostow's strong rigidity were obtained by Prasad in
the $\mathbb Q$-rank one case \cite{Pr} and
Margulis in the $\mathbb Q$-rank $\geq 2$ case \cite{Ma} 
(see also \cite{R}).  One technicality in the non-compact case is that isomorphisms 
of fundamental groups no longer induce quasi-isometries of the universal cover. 
In particular, it is no longer sufficient to just assume $\pi_1(M)\cong \pi_1(M^*)$, but
rather one needs a homotopy equivalence $f: M\rightarrow M^*$ with the property
that $f$ lifts to a quasi-isometry $\tilde f: \tM \rightarrow \tM ^*$.  We leave the
details to the interested reader.

\end{document}